# Algorithms for Groups*

## John Cannon and George Havas

School of Mathematics and Statistics
University of Sydney
NSW 2006

Key Centre for Software Technology
Department of Computer Science
University of Queensland
Queensland 4072

### Abstract

Group theory is a particularly fertile field for the design of practical algorithms. Algorithms have been developed across the various branches of the subject and they find wide application. Because of its relative maturity, computational group theory may be used to gain insight into the general structure of algebraic algorithms. This paper examines the basic ideas behind some of the more important algorithms for finitely presented groups and permutation groups, and surveys recent developments in these fields.

CR categories: I.1, G.2, G.4, F.2

# 1 INTRODUCTION

Perhaps one of the more unexpected discoveries made by mathematicians over the past two decades has been the existence of powerful algorithmic methods in many branches of algebra. Traditionally, the goal of research in algebra has been the discovery of *classification theorems* which attempt to characterize all algebraic structures satisfying a particular set of axioms. However, with the growth of interest in algebraic computation (driven, in large part, by a desire to construct symbolic solutions for various types of differential equations), mathematicians in the early 1970s were led to discover new approaches to such fundamental problems as computing the greatest common divisor of two integral polynomials and finding the irreducible factors of an integral polynomial. Since then several hundred new algorithms have been developed in various areas of algebra. It is becoming increasingly clear that most, if not all, branches of algebra have a rich algorithmic content.

Compared with most branches of algebra, the algorithmic content of group theory, *computational group theory*, has reached an advanced state of development, both in terms of the range and sophistication of the algorithms and in terms of their effectiveness in solving worthwhile problems. Thus, practical algorithms have been designed for computing detailed information concerning the structure, representations and extensions of various types of finite group. Techniques have also been developed for studying finitely-presented (fp) infinite groups. Programs implementing group theoretic algorithms find application not only in the study of groups directly but also in many

---

of the other branches of mathematics which use group theoretic methods. These include coding theory, design theory, differential equations, discrete Fourier transform theory, finite geometry, graph theory, harmonic analysis, mathematical crystallography, number theory and topology.

A group may be specified in a number of different ways: in terms of a finite presentation, as a group of permutations or matrices, or as the group of automorphisms of a combinatorial structure such as a block design, geometry or graph. Experience has shown that, as a general rule, the most powerful algorithms are those designed with a particular form of group specification in mind, *eg* permutation groups.

The major areas of activity in computational group theory are finitely presented groups, permutation groups, matrix groups, finite $p$-groups, finite soluble groups and representation theory. It is our aim to introduce the reader to some of the basic ideas that underpin the design of algorithms for group theory. A good introduction to group theory is provided by Rotman (1973). For reasons of brevity we restrict ourselves to fp-groups and permutation groups, two areas where the problems and approaches are quite different. Even within these areas, we make no attempt at completeness. A forthcoming book of Sims (to appear) provides a thorough introduction to the theory of algorithms for fp-groups, while a recent book of Butler (1991) gives an elementary introduction to computational methods for permutation groups. Both of these books include some historical information. An early account of algorithms for $p$-groups and soluble groups is given by Laue, Neubüser and Schoenwaelder (1984), while recent accounts of work in computational representation theory may be found in Lux and Pahlings (1991), Michler (1990) and Schneider (1990). Important collections of papers on computational group theory appear in Leech (1970), Atkinson (1984), Cannon (1990) and Cannon (1991b).

Computational group theory has been applied to an enormous variety of problems. Noteworthy achievements include the classification of the 4-dimensional crystallographic groups (Brown, Bülow, Neubüser, Wondratschek and Zassenhaus, 1978) and the construction of sporadic simple groups (Higman and Sims, 1968; Sims, 1973; Leon and Sims, 1977; Soicher, 1990). Many published applications prior to 1984 may be found in a bibliography by Felsch, an early version of which has been published (Felsch, 1978). A sample of more recent applications may be found in the bibliography by Cannon (1991a), which contains a partial listing of papers citing the computer algebra system Cayley.

# 2 FINITELY PRESENTED GROUPS

## 2.1 Introduction

A well-known theorem asserts that, in general, the word problem for fp-groups is undecidable (see Rotman, 1973, chapter 12). Consequently, algorithms for fp-groups are fundamentally different in nature to algorithms for groups given in some *concrete* form (*eg* permutation groups or groups of matrices over finite fields).

Let $G$ and $K$ be two fp-groups. Typical of the elementary questions mathematicians wish to answer about fp-groups are the following:

- Is $G$ the trivial group?
- Is $G$ finite?
- If $G$ is infinite, is it free?
- If $G$ is finite, what is its order and structure?
- What are the abelian (nilpotent, soluble, perfect) quotients of $G$?
- Is $G$ abelian (nilpotent, soluble, perfect)?

- Can we construct a small degree permutation representation for $G$?
- Can we construct a small degree matrix representation for $G$ over some given field?
- Are the groups $G$ and $K$ isomorphic?

The algorithms used to compute with fp-groups may be conveniently described under three headings: Todd-Coxeter or coset enumeration based methods; Knuth-Bendix or term-rewriting methods; and quotient group methods.

## 2.2   Coset Enumeration and Related Algorithms

Given an fp-group $G$, $G = \langle x_1, \ldots, x_r \mid R_1, \ldots, R_s \rangle$ (where $R_1, \ldots, R_s$ are words in the generators $x_1, \ldots, x_r$ ), and given a subgroup $H$ of $G$, $H = \langle h_1, \ldots, h_t \rangle$ (where $h_1, \ldots, h_t$ are also words in the generators), classical coset enumeration procedures attempt to construct a permutation representation for $G$, corresponding to the action of $G$ by (right) multiplication on the (right) cosets of $H$, by means of a trial-and-error process. The cosets are traditionally identified with the integers $1, \ldots, n$, where coset 1 always corresponds to the given subgroup $H$. A new coset $k$ is defined as the image of some existing coset $i$ under (right) multiplication by some generator $x_j$ of $G$ or by an inverse $x_j^{-1}$. The cosets must satisfy the following conditions:

(a) coset 1 must be mapped to itself by each of $h_1, \ldots, h_t$ ;

(b) each coset $j$ must be mapped to itself by each of the defining relators $R_1, \ldots, R_s$ and by each product $x_k x_k^{-1}$.

The action of the $G$-generators on the cosets is stored in a two-dimensional array known as a *coset table*. Enforcement of rules (a) and (b) yields values for some hitherto unknown coset table entries (*deductions*) and, also, the identification of cosets which have been multiply defined (*coincidences*). The procedure terminates when

(i) for each coset $i$, the action of each generator $x_j$ and inverse $x_j^{-1}$ on $i$ is known; and

(ii) rules (a) and (b) are satisfied.

This procedure was used extensively in hand computation prior to the development of computers. Beginning in 1952, different versions of the procedure have been adapted for machine computation and it is perhaps the most widely applied technique in computational group theory. However, despite its antiquity, our understanding of the relationship between a given presentation for $G$ and the performance of a particular version of coset enumeration when applied to that presentation is poor. Introductory descriptions of the procedure are given by Cannon, Dimino, Havas and Watson (1973), Johnson (1980), Leech (1970b, 1984) and Neubüser (1982). Sims (to appear) gives a formal account of coset enumeration in terms of automata and rational languages.

Performance of the procedure is very sensitive to changes in the rules used to introduce new cosets. Because of the many possible variations in the manner in which new cosets are introduced, coset enumeration must be regarded as a family of procedures rather than as a single procedure. For a given coset enumeration procedure, there is no computable bound, in terms of length of input and a hypothetical index, to the number of cosets which need to be defined in the coset enumeration process to complete the enumeration. (The existence of such a bound would violate the unsolvability of the word problem for finitely presented groups.) Further, Sims (to appear) has proved that there does not exist a polynomial bound, in terms of the maximum number of cosets defined, for the number of coset tables which may be constructed using simple coset table operations similar to those employed in a coset enumeration procedure. This result indicates that the running time of a coset enumeration procedure, as a function of available space, may be unpleasant.

Cannon *et al* (1973) identify a number of factors that affect the efficiency of an enumeration. The next significant study of coset enumeration procedures was carried out by Sims (to appear).

Recently, after extensive experimentation, Havas introduced new coset enumeration strategies, which exhibit dramatically better performance than previous versions when applied to many "difficult" enumerations. The performance of a coset enumeration program, in doing a particular enumeration, is measured in large part by the *maximum* number of cosets that are simultaneously defined and by the *total* number of cosets defined during the course of the enumeration. The maximum gives a direct measure of the storage requirements. For example, consider the presentation $\langle x, y, z \mid x^y x^{-3}, y^z y^{-2}, z^x z^{-4} \rangle$ for a group of order 210, essentially due to Mennicke (1959). The subgroup $\langle x \rangle$ has index 105. Methods described in Cannon *et al* (1973) define as many as a maximum of 1230970 cosets and a total of 1250191 to complete this enumeration, and the most space-economical of those methods requires a maximum of 127846 cosets and a total of 128218. The new procedure, with its default strategy, requires a maximum of 2854 cosets and total 2859, while individual tuning can lead to completion with a maximum of 1648 and total 1652. A preliminary account of this work is given by Havas (1991) and a full description is in preparation.

In favourable circumstances, the current generation of programs may successfully complete enumerations where the index of $H$ in $G$ is up to about $10^7$. Unfortunately, it is easy to construct presentations for the trivial group which will easily defeat these programs. Coset enumeration is the basis of the standard computational technique employed when attempting to prove that an fp-group $G$ is finite. If the procedure terminates, given a subgroup $H$ of $G$ known to be finite, we immediately deduce that $G$ is finite and, moreover, we obtain a bound on the order of $G$. If $G$ is not only finite but also sufficiently small so that the cosets of the trivial subgroup may be enumerated, the resulting coset table provides us with a solution to the word problem for $G$.

The range of applicability of current coset enumeration programs is limited mainly by the memory required to store the coset table. Since, in the case of non-pathological enumerations, the space required is roughly proportional to the index of $H$ in $G$, workers in the field have long dreamt of generalizing coset enumeration to a procedure capable of enumerating the double cosets $HxL$ of subgroups $H$ and $L$ of $G$. A significant step towards this goal was taken recently by Linton (1991a) who successfully implemented a double coset enumeration procedure, first suggested by Conway (1984), for the case where $H$ is a "large" subgroup of $G$ and $L$ is a small subgroup in which detailed structural computation is possible. Since the number of double cosets $HxL$ is often a small fraction of the number of single cosets $Hx$, this technique, when applicable, offers potentially great space savings.

As noted at the outset, the classical coset enumeration procedure constructs a permutation representation for $G$ on the cosets of $H$. Linton (1991b) describes a version of coset enumeration which constructs a matrix representation for $G$ over a designated field $k$ (usually a finite field). In the simplest interpretation, Linton's algorithm constructs the permutation module corresponding to the action of $G$ on the cosets of a subgroup $H$. However, since the algorithm works by constructing representations of the group algebra $kG$, given a suitable choice of ideal generators in $kG$, it is possible, for example, to directly construct constituents of a permutation module for $G$.

Our discussion so far assumes that we are given some subgroup $H$ of $G$. How do we proceed when we are unable to identify useful subgroups by direct inspection of the presentation? If $G$ is an fp-group then, for each positive integer $n$, there exist only finitely many subgroups $H$ of $G$ of index $n$. More precisely, given a homomorphism $\phi : G \to Sym(n)$, of $G$ into the symmetric group of degree $n$, such that $\phi(G)$ is transitive, then $H^\phi = \{g \in G \mid 1^{\phi(g)} = 1\}$ is a subgroup of index $n$ in $G$. Such homomorphisms may be constructed by enumerating coset tables. Given a generating set for $G$ and a positive integer $n$, there is a one-to-one correspondence between the subgroups having index $n$ in $G$ and the set of *standard* coset tables having $n$ rows (where the entries in standard tables satisfy certain ordering conditions). The *low index subgroup algorithm*, discovered independently in the sixties by Sims and Schaps, enumerates all $n$-row standard coset tables using a combination of coset enumeration and backtrack search methods (see Dietze and Schaps, 1974; Neubüser, 1982; Sims, to appear). The low index subgroup algorithm outputs either a generating set for each subgroup of index $n$, or a list of all transitive permutation representations of $G$ having degree $n$.

Versions of the low index subgroup algorithm (differing in the relative emphasis on coset enumeration or backtrack search) have been designed and implemented by Sims, Lepique (1972), Dietze and Schaps, Cannon and Gallagher (1976, unpublished), and Cannon (see Cannon and Bosma, 1991). In favourable circumstances, the low index algorithm may be used to find all subgroups having index up to 100 or even more. For example, the Cannon-Sims implementation is able to compute all conjugacy classes of subgroups in the Coxeter group $\langle a, b, c, d \mid a^2, b^2, c^2, d^2, (ab)^5, (bc)^3, (cd)^3, [a, c], [b, d], [a, d] \rangle$ with index not exceeding 240 in less than 2 hours CPU time. Since this technique is applicable to infinite groups, it provides us with a tool for proving that an fp-group is infinite. The structure of each quotient $H/H'$, where $H$ is a subgroup produced by the low index algorithm, is examined for the presence of infinite factors.

Having the coset table for a subgroup $H$ of the fp-group $G$ enables us to construct a presentation for $H$. A lemma of Schreier describes a generating set for $H$ in terms of the generators of $G$ and a system of coset representatives for $H$ in $G$. (Such coset representatives may be read off the coset table for $H$ in $G$.) The Reidemeister rewriting process permits us to rewrite the relators of $G$ and their conjugates as words in these Schreier generators. The ensuing relators constitute a presentation for $H$. Again, the coset table contains all the information needed to perform Reidemeister rewriting. Details are given by Johnson (1980) and Neubüser (1982). A variation of the Reidemeister rewriting process, which rewrites $H$-elements given as words in the generators of $G$ as words in an arbitrary generating set for $H$, was described by Benson and Mendelson (1966) (see also Neubüser, 1982). Using this theoretical basis, programs for constructing presentations of subgroups of finite index in fp-groups have been implemented by a number of workers, including Havas (1974) and Arrell and Robertson (1984).

Subgroup presentations produced by a Reidemeister-Schreier process generally involve large numbers of generators and relators and are poorly suited for human or computer use. A theorem of Teitze asserts that, given presentations for two isomorphic groups, repeated application of three simple rules (*Tietze transformations*) will suffice to transform one presentation into the other. However, there is no general algorithm for identifying the order in which the Tietze transformation rules should be applied. Presentation simplification programs which take "bad" presentations and produce "good" presentations for groups have been developed by Havas, Kenne, Richardson and Robertson (1984) and Robertson (1988).

As noted at the outset, coset enumeration over a subgroup $H$ constructs a permutation representation for $G$. We can now study the quotient of $G$ defined by this representation using permutation group methods. Conversely, given a finite group $G$ in some concrete representation, we may use the group multiplication to construct coset tables. For such a group $G$ with moderate order, Cannon (1973) shows how to construct a compact presentation for $G$ from a coset table so that we can then apply fp-group methods to $G$.

## 2.3   Term-rewriting Methods

A specialization of the Knuth-Bendix term-rewriting procedure (Knuth and Bendix, 1970) has been applied to fp-groups. Starting with a finite presentation for a monoid $M$, the Knuth-Bendix procedure for strings (*KB-procedure*) attempts to construct a *confluent* presentation for $M$. A confluent presentation for $M$ consists of a system of rewriting rules which convert any element of $M$ into a unique normal form. The KB-procedure has been studied extensively. Some applications to groups are given by Gilman (1979) and Le Chenadec (1986), while Sims (to appear) investigates the procedure in detail.

A major success of the KB-procedure in group theory was its application by Sims (1987) to verifying nilpotency of an fp-group. The nilpotent quotient algorithm (see below) is used initially to construct a polycyclic presentation for the nilpotent quotient $Q$ of $G$. Using the relations of this presentation as an initial set of rewrite rules, and a special term-ordering, Sims was able to develop an effective algorithm for verifying the triviality of the kernel of the quotient $Q$. (If the kernel is nontrivial, the algorithm fails to terminate.)

A major advantage of the KB-procedure over coset enumeration is that it may sometimes be used to construct a confluent presentation for an infinite group. In the case of a finite fp-group, coset enumeration is usually the most efficient method for constructing a confluent presentation. Epstein, Holt and Rees (1991) describe a practical algorithm based on the KB-procedure for constructing a solution to the word problem for groups known as *automatic* groups. This class of groups has solvable word problem and includes many important families of groups which arise naturally in geometry and topology (*eg* hyperbolic and Euclidean groups).

Recently, Holt and Rees (to appear a) developed a program which endeavours to determine whether or not two fp-groups $G$ and $K$ are isomorphic. The program alternates between attempting to construct an isomorphism between $G$ and $K$, and attempting to prove nonisomorphism by discovering a structural difference. The Knuth-Bendix procedure is used to construct a *reduction machine* for each of $G$ and $K$. These two reduction machines are used to systematically construct homomorphisms $\theta : G \to K$ and then to test each such homomorphism $\theta$ for the properties of being surjective and injective. The nonisomorphism testing relies on finding some structural difference by comparing various kinds of quotients of $G$ and $K$. When applied to special classes of groups, such as automatic groups, it can be quite successful. For example, Holt and Rees report that it was able to quickly settle the isomorphism question for all but two pairs in a collection of about 30 pairs of link groups. It resolved the question for the last two pairs with more difficulty, taking some hours on a Sun 3/60.

Sims (1991) employed the KB-procedure to deduce non-obvious relations in two groups. In each case, the relations could not have been discovered using the current generation of coset enumeration procedures. This is one of the first reported instances where the KB-procedure outperforms coset enumeration when both are potentially applicable.

## 2.4 Quotient Group Methods

An important technique for studying an fp-group $G$ involves constructing homomorphic images of $G$, which may be either members of some class of fp-groups having solvable word problem, permutation groups or matrix groups. As noted above, a successful coset enumeration over the cosets of some subgroup yields a homomorphism onto a permutation group, while the low index algorithm systematically searches for homomorphisms into the symmetric group $Sym(n)$, for small $n$. In this section we examine techniques for directly constructing abelian, nilpotent and soluble quotients of an fp-group. In each case, a confluent presentation for the quotient group is constructed. Note that, if the particular quotient is equal to $G$, this effectively solves the word problem for $G$.

The structure of a finitely generated abelian group $A$ may be obtained by computing the Smith normal form of its relation matrix (an integer matrix). Efficient algorithms for computing this form for large matrices have been described by Havas and Sterling (1979). If the isomorphism between $A$ and its canonical form is required, the reduced basis algorithm of Lenstra, Lenstra and Lovász (1982) may be utilized, as described by Sims (to appear, Chapter 8). Thus, given an arbitrary fp-group $G$, the structure of its maximal abelian quotient, $G/G'$, is readily obtained.

Hand computations in the sixties led to the development by Macdonald (1974) and Wamsley (1974) of algorithms for computing finite nilpotent $p$-quotients of $G$, where $p$ is a prime dividing the order of $G/G'$. Starting with an exponent-$p$-quotient of $G/G'$, the algorithms successively extend a current $p$-quotient $H$ by an elementary abelian group that is centralized by $H$. Since the extension theory is particularly simple in this situation, it is possible to design extremely effective algorithms. Nice descriptions of a basic algorithm are given by Newman (1976) and Havas and Newman (1980). The algorithm outputs the $p$-quotient in terms of a *power-commutator presentation*, a special confluent presentation. The $p$-quotient algorithm relies critically on a particular string rewriting procedure known as *commutator collection*, where the rewrite rules are the relations of a power-commutator presentation (see Felsch, 1976; Havas and Nicholson, 1976; Leedham-Green and Soicher, 1990; and Vaughan-Lee, 1990b). Since the development of the original $p$-quotient algorithms, Havas and Newman (1980) and Vaughan-Lee (1984) have introduced many improvements. A new implementation has been recently developed in Canberra by Newman and O'Brien with additional enhancements.

The algorithm has been extensively applied to the investigation of Burnside groups (see Vaughan-Lee, 1990a). As an illustration of the power of current implementations, the $p$-quotient algorithm has computed a power-commutator presentation for the three-generator restricted Burnside group of exponent 5, a group with class 17 and order $5^{2282}$. Building on the $p$-quotient algorithm, Leedham-Green and Newman designed and implemented an algorithm for generating descriptions of $p$-groups. In his PhD thesis at the Australian National University, O'Brien significantly refined these methods and successfully applied them to determine all groups with order dividing 256: there are 56092 groups of order 256 (O'Brien, 1990, 1991).

A number of variations on the original $p$-quotient algorithm have been made. Thus, a general nilpotent quotient program (with no dependency on a prime $p$ and allowing infinite quotients) was prototyped by Sims in Mathematica and a new implementation has just been developed by Nickel at ANU. Havas, Newman and Vaughan-Lee (1990) have produced an analogue of the group nilpotent quotient algorithm for graded Lie algebras. This has applications to $p$-groups where the quotients are too large to be handled by the group program and has been used to investigate questions related to the Burnside problem. Finally, during a recent visit to Australia, Vaughan-Lee developed a variation for finitely presented associative algebras.

A more difficult area is the computation of soluble quotients. Wamsley (1977), Leedham-Green (1984) and Plesken (1987) have proposed generalizations of the nilpotent quotient algorithm to a soluble quotient algorithm. The Plesken algorithm has been implemented by Wegner at St Andrews as part of his PhD thesis and has had some success. Baumslag, Cannonito and Miller (1981) outlined a method for constructing polycyclic quotients. While their interest was purely theoretical, Sims (1990b) has further developed their ideas and implemented them in the special case of metabelian quotients. An implementation of the general algorithm (see Sims, to appear, chapters 9 and 10) involves a sophisticated combination of many algorithms, including Gröbner basis techniques (Buchberger, 1985). Success has been reported in specific cases by Neubüser and Sidki (1988), Newman and O'Brien (to appear) and Havas and Robertson (to appear), all using combinations of previously described programs.

At the DIMACS Workshop on Groups and Computation in 1991, Holt and Rees (to appear b) demonstrated a program for finding certain quotients of finitely presented groups. A backtrack search attempts to construct a homomorphism between the given fp-group $G$ and nominated permutation representations of selected finite groups. In particular, this program may be used to identify small perfect groups that occur as quotients of $G$. This work builds on the classification by Holt and Plesken (1989) of all perfect groups of order up to a million. Having found a representation of $G$, the program converts it to a regular representation, and then attempts to construct larger quotients of $G$ by performing elementary abelian extensions using either a Reidemeister-Schreier or $p$-quotient algorithm. The portion of the subgroup lattice obtained is represented graphically and may be manipulated interactively.

# 3    PERMUTATION GROUPS

## 3.1    Introduction

Computational methods for permutation groups differ fundamentally from those designed for fp-groups since the undecidability of the word problem means that many fp-group procedures cannot be guaranteed even to terminate. Thus, if we are lucky, the fp-procedures may provide us with some information about the overall properties of the group under consideration. However, in the case of a permutation group, it is possible to design fast algorithms which can provide us with highly detailed information about the structure of the group.

Over the past two decades, a large number of useful permutation group algorithms have been discovered. Because the field is so extensive, this survey is restricted to a few fundamental classes of algorithms and no attempt is made at completeness, even in the case of the areas we do consider. A more complete account of basic permutation group algorithms is given by Cannon (1984) and Bosma and Cannon (to appear). Butler (1991) provides an introductory account of many basic algorithms. We do not attempt to survey the extensive literature on the complexity of permutation group algorithms but rather refer the reader to Babai (1991).

Let $G$ be a permutation group acting on the finite set $\Omega$ and suppose $G$ is given in terms of a small set $X$ of generating permutations. The following are representative of the type of information sought by permutation group theorists:

• What are the $G$-invariant subsets (partitions) of the set $\Omega$?

• What is the order of the group $G$?

• Is a given element of $Sym(\Omega)$ a member of the group $G$?

• Find generators for the stabilizer of a sequence (set) of elements of $\Omega$.

• Determine the various series of characteristic subgroups of $G$: derived series; lower central series; upper central series.

• Find generators for the Sylow $p$-subgroup of $G$, where $p$ is a prime dividing the order of $G$.

• Compute the centralizer (normalizer, normal closure) of a subgroup $H$ of $G$.

• Determine representatives for the conjugacy classes of elements of $G$.

• Compute a composition series for $G$ and determine the isomorphism type of each composition factor.

The construction of a $G$-invariant partition of $\Omega$ proceeds by computing the finest $G$-invariant partition that contains a given pair of points $\alpha$ and $\beta$. Consideration of all such pairs of points in turn will yield the complete list of minimal $G$-invariant partitions of $\Omega$. Various reductions which greatly improve efficiency are possible. Thus, if $G$ is transitive, it suffices to consider the pairs $\alpha$, $\beta$, where $\alpha$ is fixed and $\beta$ runs through the elements of the set $\Omega - \{\alpha\}$. A particularly efficient version of such an algorithm is described by Atkinson (1975). Used to test primitivity of $G$, the algorithm has running time $O(mn^2)$, where $m = |X|$ and $n = |\Omega|$.

## 3.2  Base and Strong Generating Sets

All but the first of the above questions involve quantifying over the *set* of elements of $G$. The design of effective algorithms depends critically on the method chosen to represent this set. The representation should have the following features:

(a) it should display key aspects of the group structure;

(b) it should have the property of inheritance, thus a subgroup of $G$ should directly inherit its representation from that of $G$;

(c) it should be effectively computable.

Sims (1970, 1971a, 1971b) introduced the notion of a base and strong generating set (*BSGS*) as the appropriate set representation for a permutation group. A *base* for $G$ is a sequence $B = \langle \beta_1, \ldots, \beta_k \rangle$ of distinct elements of $\Omega$ such that the identity is the only element of $G$ that fixes $B$ pointwise. Let $G_{\alpha_1, \ldots, \alpha_k}$ denote the pointwise stabilizer of $\{\alpha_1, \ldots, \alpha_k\} \subset \Omega$. Then $B$ defines a sequence of subgroups, $G = G^{(1)} > G^{(2)} > \ldots > G^{(k+1)} = \langle 1 \rangle$, where $G^{(i)} = G_{\beta_1, \ldots, \beta_{i-1}}$. A *strong generating set* $S$ for $G$, relative to the base $B$, is a generating set for $G$ which contains generators for each subgroup in the chain. Given $B$ and $S$, it is a straightforward matter to compute the orbit $\Delta_i = \beta_i^{G^{(i)}}$ and a transversal $U^{(i)}$ for $G^{(i+1)}$ in $G^{(i)}$, for $i = 1, \ldots, k$.

Knowledge of $\Delta_i$ and $U^{(i)}$ for $i = 1, \ldots, k$ immediately gives us the order of $G$, a membership test for $G$, and the possibility of listing the elements of $G$ without repetition. The central philosophy in permutation group computation is to assume that the group $G$ is represented in terms of a BSGS (whose construction is discussed below). Sims showed that the availability of a BSGS enables us to design effective algorithms for constructing important classes of subgroups and quotient groups of $G$. Wherever possible, algorithms for constructing such subgroups and quotient groups of $G$ are organized in such a way that the new group *inherits* a BSGS from the known BSGS for $G$. For example, the backtrack algorithms for constructing subgroups, described below, automatically produce a BSGS for the subgroup.

The design of fast methods for constructing a BSGS for $G$ is one of the first problems which must be solved. In 1967, Sims developed a BSGS algorithm based on the following lemma of Schreier:

Lemma: Let $G$ be a group with generating set $X$, and let $H$ be a subgroup of $G$. If $U$ is a right transversal for $H$ in $G$ such that $U$ contains the identity element then $H$ is generated by the set $\{ux\phi(ux)^{-1} \mid u \in U, x \in X\}$, where $\phi(ux)$ is the unique element of $U$ such that $Hux = H\phi(ux)$.

In theoretical terms, the algorithm proceeds as follows: For $\beta_1$ choose any point not fixed by every element of $X$. Compute the orbit of $\beta_1$ together with the corresponding transversal $U^{(1)}$. Now apply Schreier's lemma to construct a set of generators for $G_{\beta_1} = G^{(2)}$. Iterating this process, we successively construct base points $\beta_2, \ldots, \beta_k$, and generating sets for the stabilizers $G^{(3)}, \ldots, G^{(k+1)}$. The process terminates when we reach the trivial subgroup. This algorithm is impractical as it stands because of the large number of generators produced by the lemma. In fact, the number of elements needed to generate the stabilizer is usually a tiny fraction of the number produced by the lemma. (It is easy to see that at most $n$ elements are needed to generate the stabilizer of a point in $G$.) With careful organization, the above idea yields a BSGS algorithm that produces small generating sets for the terms of the stabilizer chain. In practice, the algorithm works well for degrees less than 100 but becomes impractical as the degree increases beyond a few hundred. For the case of a soluble permutation group, Sims (1990a) describes a BSGS algorithm which takes account of the particular structure of these groups.

For larger degrees a different approach is adopted. First a "probable" BSGS for $G$ is constructed, and then an algorithm is applied which either *verifies* that the BSGS is correct, or establishes that it is incomplete. A "probable" BSGS may be constructed very quickly by using a fixed number of randomly chosen elements of $G$ in place of the products $ux$ in Schreier's lemma (*probabilistic Schreier algorithm*). The main inductive step in BSGS verification is the following: Suppose $H$ is a subgroup of $G^{(2)}$ with a certified BSGS. If we can show that $H = G^{(2)}$, then we will have verified the correctness of the BSGS for $G$. In practice, $H$ will be the approximation to $G^{(2)}$ constructed by the probabilistic Schreier algorithm. The first verification algorithm was developed by Sims and published by Leon (1980) and involves using the Todd-Coxeter algorithm to construct a presentation for $G$ in terms of strong generators. A presentation is assumed to be known for $H$ on its strong generating set by induction. The verification involves comparing the index of $H$ in $G$ with the length of the orbit $\beta_1^G$. This method made it practical to construct BSGSs for groups having degree in the low thousands. In 1986, Brownie and Cannon implemented a new verification algorithm suggested by Sims. Instead of testing each of the Schreier generators given in the above lemma for membership of $H$, we test a much smaller subset defined in terms of representatives for the orbits of certain one-point and two-point stabilizers. This algorithm, as implemented in Cayley, has constructed BSGSs for groups of degree up to 500000 and, given sufficient memory, is capable of computing BSGSs for groups having degrees up to a million.

The theoretical complexity of BSGS algorithms depends heavily upon the choice of data structures, particularly on the data structure used to store the transversals $U^{(i)}$. Sims stored each transversal $U^{(i)}$ in the form of a linearized tree structure known as a *Schreier vector* and the running time of his original algorithm was bounded by $O(n^6)$, where $n$ is the degree of $G$. Using the so-called *labelled branching* data structure for the $U^{(i)}$, Jerrum (1986) described a variant of the original Sims algorithm with running time $O(n^5)$. Cooperman and Finkelstein (1991) describe an algorithm which verifies strong generation in $O(n^4)$ time. More experimental work needs to be done in order to establish whether or not the Cooperman-Finkelstein algorithm is a practical competitor to the Brownie-Cannon-Sims algorithm.

Given a particular base $B$ for $G$ and strong generators relative to $B$, there exist fast algorithms for computing strong generators for $G$ relative to a different base $B'$ for $G$. Sims (1971a, 1971b) shows how to modify the strong generators when two adjacent base points are interchanged. By concatenating the new base $B'$ onto the end of $B$, the desired base change may be achieved by performing a succession of adjacent interchanges. Some modifications to speed up this procedure are described in Butler (1991). Cooperman, Finkelstein and Luks (1989) describe a base change algorithm with running time $O(kn^2)$. A different approach to changing base is to use a probabilistic procedure similar to the random Schreier algorithm mentioned above (see Leon, 1980). Such a procedure is employed in the Brownie-Cannon-Sims verification algorithm. The complexity of such algorithms has been analyzed by Cooperman, Finkelstein and Sarawagi (1990).

## 3.3 Subgroup Constructions

The base change algorithm enables us to rapidly construct a BSGS for the pointwise stabilizer of a sequence of points from $\Omega$. The availability of a BSGS for $G$ enables us to test membership of $G$ in polynomial time, where the degree of the polynomial (dependent on the data structure and space used) may be as low as 2. The combination of a BSGS algorithm and membership testing immediately allows us to construct the normal closure of a subgroup of $G$. The availability of a normal closure algorithm in turn enables us to construct the derived series and lower central series for $G$, and also allows us to test subgroups for the properties of being perfect, soluble, nilpotent, subnormal, etc. (see Butler and Cannon, 1982).

The availability of the BSGS representation of a permutation group provides the appropriate foundation for the design of efficient backtrack searches for subgroups of $G$ whose elements satisfy some elementary property. By carefully choosing a base adapted to a particular backtrack search, we can often greatly reduce the size of the search tree. The BSGS backtrack search of a permutation group was introduced by Sims in 1970, when he described backtrack algorithms for computing centralizers and intersections of subgroups. Over the next decade, Butler, Cannon and Sims developed backtrack searches for constructing set stabilizers, normalizers, Sylow $p$-subgroups and for testing conjugacy of elements and subgroups (see Sims, 1971a, 1971b; Butler, 1982, 1983, 1991). Holt (1991) presents a backtrack algorithm for computing subgroup normalizers which employs many additional tests to prune the backtrack search tree. The performance of his algorithm is superior in many cases to the Butler algorithm as implemented in the Cayley library of intrinsic functions.

Very recent work of Leon (1991) represents a major step in the evolution of backtrack algorithms for permutation groups. The efficiency of a backtrack search is heavily dependent upon the information available to prune the search tree. Using the idea of successive refinement of ordered partitions, first introduced by McKay (1978, 1981) as part of his highly successful graph isomorphism algorithm, Leon is able to devise new and powerful tests. An early implementation of a set stabilizer algorithm based on these ideas demonstrates performance that is dramatically superior to the "first generation" set stabilizer algorithm. As a result of this work we can expect a new generation of backtrack algorithms, exhibiting superior performance, to emerge in the near future.

Let $p$ be a prime dividing $|G|$, and let $P$ denote the Sylow $p$-subgroup of $G$. Traditional approaches to computing $P$ have involved performing a series of cyclic extensions commencing with the subgroup generated by an element of $p$-power order. Butler and Cannon (1989) describe an algorithm which employs a backtrack search for possible extending elements. This algorithm is restricted to groups having degrees in the low hundreds and primes whose exponent in $|G|$ does not exceed 16. Holt (1991) suggests using his fast normalizer algorithm to locate possible extending elements. Butler and Cannon (1991) present a recursive method based on reduction of the degree. The reduction is based on the observation that if $z$ is a $p$-central element having order $p$, then the kernel of the action of the centralizer of $z$ on the cycles of $z$ is a $p$-group. This algorithm, which may involve a probabilistic search for a $p$-central element of $G$, generally runs a great deal faster than the backtrack method. Yet another method based on the use of homomorphisms has been suggested by Atkinson and Neumann (1990). Kantor (1985) and Kantor and Taylor (1988) give polynomial time algorithms for computing Sylow $p$-subgroups which appear to be of theoretical interest only.

## 3.4   Abstract Structure

Given the availability of constructions for centralizer, normal closure and Sylow subgroup, one can contemplate computing a description of the abstract structure of $G$. Such a goal became particularly attractive with the completion of the classification of the finite simple groups.

A very important requirement is to be able to determine quickly whether or not a permutation group $G$ of degree $n$ contains the alternating group $Alt(n)$ in its natural representation. Indeed, it is desirable to recognize this situation before incurring the expense of constructing a BSGS for $G$. The starting point is a theorem of Jordan which states that if a primitive group contains a $p$-cycle, where $p \leq n - 3$ is a prime, then $G$ contains $Alt(n)$. This theorem has been subsequently extended by other workers to identify many other permutation cycle types that are *Jordan witnesses* for $Alt(n)$. Cameron and Cannon, in work extending over a decade, have constructed a subtle probabilistic recognition procedure for groups containing $Alt(n)$. The basic idea is to sample a very small number of elements of $G$ and, if $G$ does indeed contain $Alt(n)$, deduce from these elements the primitivity of $G$ and the presence of Jordan witnesses. This algorithm is applicable to groups having degree up to a million and, if $G$ does contain $Alt(n)$, it usually only needs to sample two or three elements in order to recognize the fact.

As a next step, Cameron and Cannon (1991) developed an algorithm for identifying any doubly transitive group. By carefully analyzing the lengths of the orbits of a two-point stabilizer (three-point stabilizer in a triply transitive group), the algorithm avoids having to compute the derived subgroup of $G$, except in the case of some relatively small groups (*eg* one-dimensional affine groups).

At the 1985 Groups — St Andrews meeting, Neumann (1986) described a practical algorithm, based on the O'Nan-Scott Theorem, for determining a BSGS for each composition factor of a general permutation group $G$. The general strategy involves reducing to a primitive group $T$, locating the socle of $T$, and then splitting the socle into its simple direct factors. Luks (1987) published a polynomial time algorithm for this problem which appears to be of theoretical interest only. Neumann's algorithm made some use of the fact that practical computation with such an algorithm will be restricted to groups having degree at most a few million. Kantor (1991) pushes this approach somewhat further and shows that, except in the case of a few small groups, it is possible to name the composition factors of any permutation group having degree not exceeding $10^6$ at the cost of computing a BSGS for $G$ and the derived subgroup of $G$. Thus, Kantor improves on the Neumann method by avoiding the step of constructing and splitting the socle of $T$. It appears that the Kantor algorithm may be generalized to groups of degree up to $10^8$ so that his approach covers all permutation groups that are presently amenable to practical computation.

Going in the other direction, Holt has developed an algorithm for constructing the first and second cohomology groups of a permutation group. This very complex algorithm builds on several of the permutation group algorithms described above, employs the nilpotent quotient algorithm and involves intricate module-theoretic calculations (see Holt, 1984, 1985a, 1985b).

# 4 ACKNOWLEDGEMENTS

Both authors were partially supported by Australian Research Council grants. The second author was partially supported by DIMACS/Rutgers-NSFSTC88-09648.

## BIOGRAPHICAL NOTES

John Cannon is a Reader in Pure Mathematics at the University of Sydney. He is leader of the Computer Algebra Group there and the principal architect of the computer algebra system Cayley.

George Havas is a Senior Lecturer in Computer Science at the University of Queensland. His previous appointments include: Lecturer in Computing at the Canberra College of Advanced Education (now University of Canberra); Research Fellow in Mathematics at the Institute of Advanced Studies, Australian National University; Principal Research Scientist, Division of Computing Research, CSIRO; and Director of Research and Development, Csironet.